\documentclass{amsart}
\usepackage{amssymb,amsfonts,amsmath}
\usepackage{graphics,color}
\usepackage[matrix,arrow]{xy}

\numberwithin{equation}{section}

\theoremstyle{plain}

\newtheorem{thm}{Theorem}[section]
\newtheorem{theorem}[thm]{Theorem}

\newtheorem{corollary}[thm]{Corollary}
\newtheorem{lemma}[thm]{Lemma}

\newtheorem{proposition}[thm]{Proposition}

\newtheorem{remark}[thm]{Remark}

\newtheorem{example}[thm]{Example}

\newcommand{\R}{\mathbb{R}}

\newcommand{\C}{\mathbb{C}}

\newcommand{\htp}{\simeq}
\newcommand{\iso}{\cong}

\newcommand{\A}{\mathcal{A}}
\newcommand{\B}{\mathcal{B}}
\newcommand{\CC}{\mathcal{C}}
\newcommand{\D}{\mathcal{D}}
\newcommand{\PP}{\mathcal{P}}

\newcommand{\K}{\mathbb{K}}

\newcommand{\Coh}{\mathfrak{Coh}}
\newcommand{\Fuk}{\mathfrak{Fuk}}

\begin{document}
\title[Lefschetz fibrations]{Suspending Lefschetz fibrations, \\ with an application to Local Mirror Symmetry}
\author{Paul Seidel}
\date{July 12, 2009}
\maketitle

\section{Introduction}

Let $Y$ be a smooth toric del Pezzo surface, and $K_Y$ the total space of its canonical bundle. Let $D^b(\Coh(K_Y))$ be the bounded derived category of coherent sheaves on $K_Y$, and $D^b_Y(\Coh(K_Y))$ the full subcategory consisting of complexes whose cohomology is supported on the zero-section $Y \subset K_Y$. The mirror to $K_Y$ (see \cite[Section 8]{hori-iqbal-vafa00} or \cite{gross01b}) is the hypersurface
\begin{equation}
H = \{(x,y) \in (\C^*)^2 \times \C^2 \;:\; y_1y_2 + p(x) = z\},
\end{equation}
Here, $p: (\C^*)^2 \rightarrow \C$ is the superpotential mirror to $Y$ (following \cite{givental94} or \cite{hori-vafa00}), and $z$ is any regular value of $p$. $H$ is an affine threefold with trivial canonical bundle. Hence, it has a Fukaya category $\Fuk(H)$, whose objects are compact exact Lagrangian submanifolds equipped with gradings and $\mathit{Spin}$ structures. This is an $A_\infty$-category over $\C$. Consider the associated derived category $D(\Fuk(H))$.

\begin{theorem} \label{th:local-hms}
There is a full embedding of triangulated categories,
\begin{equation} \label{eq:emb}
D^b_Y(\Coh(K_Y)) \hookrightarrow D(\Fuk(H)).
\end{equation}
\end{theorem}

To explain the strategy of the proof, we need to return to $Y$ itself. Homological mirror symmetry for such varieties, first considered by Kontsevich \cite{kontsevich98}, is now a theorem (see \cite{seidel00b} for $Y = \C P^{2}$, \cite{auroux-katzarkov-orlov04} for several more cases, and \cite{auroux-katzarkov-orlov05,ueda06} for all remaining ones). On the other hand, Segal \cite{segal07} and Ballard \cite{ballard08} independently proved a result which describes, on the level of the underlying $A_\infty$-algebras, the relation between $D^b(\Coh(Y))$ and $D^b_Y(\Coh(K_Y))$. The missing piece, which we provide here, is a corresponding statement about the relation between the Fukaya category of the Lefschetz fibration $p$ and that of the hypersurface $H$. We can prove this relation only partially, which is why we get an embedding \eqref{eq:emb} instead of the ultimately desirable equivalence.

The proof of Theorem \ref{th:local-hms} is in fact an application of a more general result about suspending Lefschetz fibrations (Theorem \ref{th:suspension}), which means passing from $p(x)$ to $p(x) - y_1^2 - \cdots - y_m^2$ for any number of variables. Independently of the present work, Futaki and Ueda \cite{futaki-ueda09} have applied similar methods to $p(x) - y^k$ for $k \geq 2$. However, their interest is directed towards mirror symmetry for Landau-Ginzburg models, hence their computations concern only directed Fukaya categories.

{\em Acknowledgments.} This problem was brought to my attention by Matthew Ballard. I would like to thank him, Mohammed Abouzaid, Mark Gross, and Ivan Smith for useful conversations. NSF grant DMS-0652620 provided partial support.

\section{Algebraic suspension\label{sec:def}}

Fix a ground field $\K$, which will be used throughout. We will describe a suspension construction, which starts from a pair $(\A,\B)$ consisting of an $A_\infty$-algebra $\B$ and an $A_\infty$-subalgebra $\A \subset \B$, and produces another pair $(\A^\sigma,\B^\sigma)$ of the same kind. In fact $\A^\sigma$ is isomorphic to $\A$ in a straightforward way, but $\B^\sigma$ is usually not even quasi-isomorphic to $\B$. We will see that in general, the construction leads to a loss of information, which means that it cannot be reversed in a meaningful (that is, quasi-isomorphism invariant) way.

Here is the most elementary description. As a graded vector space
\begin{equation}
\B^\sigma = \A_+ \oplus \A_- \oplus \B[-1],
\end{equation}
where both $\A_+$ and $\A_-$ are copies of $\A$. The shift $[-1]$ means that elements of $\B^\sigma$ of degree $r$ are triples $(a_+,a_-,b)$, where $a_\pm \in \A$ have degree $r$, and $b \in \B$ has degree $r-1$. The differential is
\begin{equation} \label{eq:mu1-stabilized}
\mu^1_{\B^\sigma}(a_+,a_-,b) = (\mu_\A^1(a_+),\mu_\A^1(a_-), -\mu_\B^1(b) - a_+ + a_-).
\end{equation}
The higher order $A_\infty$-structure maps, for $d \geq 2$, are
\begin{multline} \label{eq:mu2-stabilized}
\mu^d_{\B^\sigma}((a_{d,+},a_{d,-},b_d),\dots,(a_{1,+},a_{1,-},b_1)) \\ =
\Big(\mu^d_\A(a_{d,+},\dots,a_{1,+}),\mu^d_\A(a_{d,-},\dots,a_{1,-}), \\
\sum_{i=1}^d (-1)^{\|a_{1,-}\| + \cdots + \|a_{i-1,-}\|+1} \mu_\B^d(a_{d,+},\dots,a_{i+1,+},b_i,a_{i-1,-},\dots,a_{1,-})\Big).
\end{multline}
The notational convention is that $\|a\| = \mathrm{deg}(a) - 1$ is the reduced degree. The subspace $\A^\sigma \subset \B^\sigma$ just consists of triples of the form $(a,a,0)$.

\begin{example} \label{th:dga}
It maybe makes sense to write down the construction in the simpler special case of dga's, taking into account the differences in sign conventions. Namely, if $\A \subset \B$ are dga's, the dga structure on the suspension $\B^\sigma$ is given by
\begin{equation}
\begin{aligned}
& d(a_+,a_-,b) = (da_+,da_-,db - (-1)^{\deg(a_+)} a_+ + (-1)^{\deg(a_-)} a_-), \\
& (a_{2,+},a_{2,-},b_2) \cdot (a_{1,+},a_{1,-},b_1) = \\ & \qquad = (a_{2,+}a_{1,+},a_{2,-}a_{1,-},a_{2,+}b_1 + (-1)^{\mathrm{deg}(a_{1,-})} b_2 a_{1,-}).
\end{aligned}
\end{equation}
\end{example}

The fact that \eqref{eq:mu1-stabilized}, \eqref{eq:mu2-stabilized} satisfy the $A_\infty$-associativity equations can be checked by hand, but there is also an equivalent description which makes these equations obvious. Start with the following contractible complex of vector spaces:
\begin{equation}
C = \{\cdots 0 \rightarrow \K \xrightarrow{1} \K \rightarrow 0 \cdots\},
\end{equation}
where the nontrivial generators are in degrees $-1$ and $0$. Let $\mathit{hom}_\K(C,C)$ be the differential graded algebra of endomorphisms of $C$. Consider the tensor product $\B \otimes \mathit{hom}_\K(C,C)$, which is again an $A_\infty$-algebra.
%(for the sign conventions see \cite[Remark 1.11]{seidel04}).
As a graded vector space,
\begin{equation}
\B \otimes \mathit{hom}_\K(C,C) = \B[1] \oplus \B_+ \oplus \B_- \oplus \B[-1],
\end{equation}
where $\B_+$ is $\B$ tensored with the endomorphisms of the degree $0$ part of $C$, and $\B_-$ similarly for the degree $-1$ part of $C$. Define $\B^\sigma$ as an $A_\infty$-subalgebra of $\B \otimes \mathit{hom}_\K(C,C)$ in the obvious way. This leads to the formulae described above. For instance, the additional terms in the differential \eqref{eq:mu1-stabilized} are precisely the ones inherited from the differential on $\mathit{hom}_\K(C,C)$.

\begin{remark} \label{th:twisted}
One can recast the previous description in slightly different language as follows. Assume that $\B$ is strictly unital, with identity element $e$. Consider it as an $A_\infty$-category with a single object $V$, and let $\mathit{Tw}(\B)$ be the associated $A_\infty$-category of twisted complexes. This contains a contractible object, which is the mapping cone $S = \mathit{Cone}(e: V \rightarrow V)$. Then
\begin{equation}
\mathit{hom}_{\mathit{Tw}(\B)}(S,S) = \B \otimes \mathit{hom}_\K(C,C)
\end{equation}
as an $A_\infty$-algebra, and one can again define $\B^\sigma$ as an $A_\infty$-subalgebra of this (our sign conventions differ from those of \cite[Section 3k]{seidel04}, essentially because we think of twisted complexes as formal tensor products with graded vector spaces on the right side, rather than the left one).
%We refer to \cite[Section 3k]{seidel04} for the sign conventions used in setting up twisted complexes, and to %\cite[Section 3p]{seidel04} for the definition of the mapping cone.
\end{remark}

There is a special class of examples where the effect of suspension is quite simple. To explain this, we need the notion of an $A_\infty$-bimodule over $\A$, which is a graded vector space $\PP$ together with operations
\begin{equation} \label{eq:mu-p}
\mu_\PP^{s|1|r}: \A^{\otimes s} \otimes \PP \otimes \A^{\otimes r} \longrightarrow \PP[1-r-s]
\end{equation}
for all $r,s \geq 0$, satisfying appropriate equations. A straightforward example is the diagonal bimodule, which is $\A$ itself with
\begin{equation} \label{eq:diagonal}
\mu_{\A}^{d-i|1|i-1}(a_d,\dots,a_1) \\ = (-1)^{\|a_1\| + \cdots + \|a_{i-1}\|+1}
\mu_\A^d(a_d,\dots,a_1).
\end{equation}
We will also need the shift operation, which takes $\PP$ to the shifted vector space $\PP[-1]$ with operations
\begin{multline}
\mu_{\PP[-1]}^{s|1|r}(a_{r+s},\dots,a_{r+1},p,a_r,\dots,a_1) \\ = (-1)^{\|a_1\| + \cdots + \|a_r\|+1}
\mu_{\PP}^{s|1|r}(a_{r+s},\dots,a_{r+1},p,a_r,\dots,a_1).
\end{multline}
Hence, the shifted diagonal bimodule $\A[-1]$ has $\mu_{\A[-1]}^{s|1|r} = \mu_\A^{s+1+r}$. For a more thorough exposition of the theory of $A_\infty$-bimodules and their morphisms, see \cite{tradler01,seidel08} (our sign conventions are related to those in \cite{seidel08} by reversing the ordering of the entries in \eqref{eq:mu-p}, which brings things in line with \cite{seidel04}).

For any given $A_\infty$-algebra $\A$ and $\A$-bimodule $\PP$, one can define the trivial extension $A_\infty$-algebra $\A \oplus \PP$, which is the direct sum of vector spaces equipped with the structure maps
\begin{multline} \label{eq:trivial-extension}
\mu^d_{\A \oplus \PP}((a_d,p_d),\dots,(a_1,p_1)) = \Big(\mu_\A^d(a_d,\dots,a_1), \\
\sum_{i=1}^d (-1)^{\|a_1\| + \cdots + \|a_{i-1}\|+1} \mu_\PP^{d-i|1|i-1}(a_d,\dots,a_{i+1},p_i,a_{i-1},\dots,a_1)\Big).
\end{multline}
This obviously contains $\A$ as an $A_\infty$-subalgebra. The suspension of $\B = \A \oplus \PP$ is $\B^\sigma = \A_+ \oplus \A_- \oplus \A[-1] \oplus \PP[-1]$, with
\begin{equation}
\mu^1_{\B^\sigma}(a_+, a_-,a,p) = (\mu_\A^1(a_+),\mu_\A^1(a_-), -\mu_\A^1(a) - a_+ + a_-, \mu_\PP^{0|1|0}(p))
\end{equation}
and, for $d \geq 2$,
\begin{multline}
\mu^d_{\B^\sigma}((a_{d,+},a_{d,-},a_d,p_d),\dots,(a_{1,+},a_{1,-},a_1,p_1)) \\ =
\Big(\mu^d_\A(a_{d,+},\dots,a_{1,+}),\mu^d_\A(a_{d,-},\dots,a_{1,-}), \\
\sum_{i=1}^d (-1)^{\|a_{1,-}\| + \cdots + \|a_{i-1,-}\|+1} \mu_\A^d(a_{d,+},\dots,a_{i+1,+},a_i,a_{i-1,-},\dots,a_{1,-}), \\
\sum_{i=1}^d \mu_\PP^{d-i|1|i-1}(a_{d,+},\dots,a_{i+1,+},p_i,a_{i-1,-},\dots,a_{1,-})
\Big).
\end{multline}
This contains the $A_\infty$-subalgebra of elements of the form $(a,a,0,p)$, which is isomorphic to the trivial extension $\A \oplus \PP[-1]$. Moreover, the induced differential on the quotient $\B^\sigma/(\A \oplus \PP[-1]) \iso \A \oplus \A[-1]$ is acyclic. We have therefore shown:

\begin{lemma} \label{th:trivial-extension}
If $\B$ is a trivial extension $\A \oplus \PP$, then its suspension $\B^\sigma$ is quasi-isomorphic to the trivial extension $\A \oplus \PP[-1]$. \qed
\end{lemma}

\section{A topological viewpoint\label{sec:topology}}

This section is a digression from our main argument, and its only purpose is to build some topological intuition.
Let $U$ be a topological manifold with boundary $W = \partial U$ (in fact, more general topological spaces and subspaces are also possible). Let $\A = C^*(U)$ be the dga of singular cochains with coefficients in $\K$. Set $\B = C^*(U) \oplus C^*(U,W)[1]$, with the differential and multiplication given, in terms of the standard differential $\delta$ and cup-product on cochains, by
\begin{equation}
\begin{aligned}
& d(b,c) = (\delta b + (-1)^{\mathrm{deg}(c)} c, \delta c), \\
& (b_2,c_2) \cdot (b_1,c_1) = (b_1 b_2,b_2 c_1 + (-1)^{\mathrm{deg}(b_1)} c_2 b_1).
\end{aligned}
\end{equation}
By definition, the map $\B \rightarrow C^*(W)$, $(b,c) \mapsto b|W$, is a quasi-isomorphism of dga's. The point of replacing $C^*(W)$ by its quasi-isomorphic version $\B$ is that it allows us to view the restriction map $C^*(U) \rightarrow C^*(W)$ as the inclusion $\A \hookrightarrow \B$.

Consider the suspension of $\B$, which following Example \ref{th:dga} is the dga $\B^\sigma = C^*(U) \oplus C^*(U) \oplus C^*(U)[-1] \oplus C^*(U,W)$ with
\begin{equation}
\begin{aligned}
& d(a_+,a_-,b,c) = (\delta a_+,\delta a_-, \\ & \delta b - (-1)^{\mathrm{deg}(a_+)} a_+ + (-1)^{\mathrm{deg}(a_-)} a_- + (-1)^{\mathrm{deg}(c)} c, \delta c), \\
& (a_{2,+},a_{2,-},b_2,c_2) \cdot (a_{1,+},a_{1,-},b_1,c_1) = \\ & =
(a_{2,+}a_{1,+},a_{2,-}a_{1,-}, a_{2,+}b_1 + (-1)^{\mathrm{deg}(a_{1,-})} b_2 a_{1,-}, a_{2,+} c_1 + c_2 a_{2,-}).
\end{aligned}
\end{equation}
Let $W^\sigma = U_+ \cup_W U_-$ be the manifold obtained by gluing two copies of $U$ along their common boundary. This comes with a dga homomorphism
\begin{equation} \label{eq:sandwich}
\begin{aligned}
& C^*(W^\sigma) \longrightarrow \B^\sigma, \\
& a \longmapsto (a|U_+,a|U_-,0,a|U_+ - a|U_-).
\end{aligned}
\end{equation}
This implicitly uses the fact that the restrictions $a|U_\pm$ agree along their boundaries, so that their difference is naturally a relative cochain. As a map of chain complexes, \eqref{eq:sandwich} sits in a commutative diagram with short exact rows,
\begin{equation}
\xymatrix{
 0 \ar[r] & C^*(W^\sigma,U_-) \ar[r] \ar[d] & C^*(W^\sigma) \ar[d] \ar[r] & C^*(U_-) \ar[d] \ar[r] & 0 \\
 0 \ar[r] & C^*(U,W) \ar[r] & \B^\sigma \ar[r] & C^*(U) \oplus C^*(U) \oplus C^*(U)[-1] \ar[r] & 0.
}
\end{equation}
The lower row consists of the maps $c \mapsto (c,0,0,c)$ and $(a_+,a_-,b,c) \mapsto (a_+ - c,a_-,b)$. The left hand vertical arrow is restriction from $W^\sigma$ to $U_+ = U$, and the right hand one is $a \mapsto (a,a,0)$ with respect to the identification $U_- = U$. The two last-mentioned maps are quasi-isomorphisms, hence so is \eqref{eq:sandwich}.

Next, $W^\sigma$ is naturally the boundary of a manifold $U^\sigma$ which is homotopy equivalent to $U$, namely
\begin{equation}
U^\sigma = U \times [-1,1] \; / \; \text{$(x,t) \sim (x,s)$ for $x \in W$}.
\end{equation}
The maps induced by projection $U^\sigma \rightarrow U$ and inclusion $W^\sigma \hookrightarrow U$ fit into a commutative diagram of dga's,
\begin{equation}
\xymatrix{
\ar[d]^{\htp} C^*(U) \ar[r]^-{\htp} & C^*(U^\sigma) \ar[r] & C^*(W^\sigma) \ar[d]^{\htp} \\
\A^\sigma \ar[rr] && \B^\sigma.
}
\end{equation}
Hence, in this context suspension is the algebraic counterpart of the topological process of passing from $(U,W)$ to $(U^\sigma,W^\sigma)$.

\section{General properties\label{sec:ext}}

We return to general pairs $(\A,\B)$. The inclusion $\A \subset \B$ makes $\B$ into an $A_\infty$-bimodule over $\A$, in the same way as in \eqref{eq:diagonal}. It contains the diagonal bimodule $\A$, so the quotient $\B/\A$ is naturally an $A_\infty$-bimodule as well, and we have a short exact sequence
\begin{equation}
0 \rightarrow \A \xrightarrow{\iota} \B \xrightarrow{\pi} \B/\A \rightarrow 0.
\end{equation}

\begin{lemma} \label{th:stabilized-natural}
Suppose that we have two pairs $(\A,\B)$ and $(\A,\tilde\B)$, together with an $\A$-bimodule quasi-isomorphism $\phi: {\tilde\B} \rightarrow \B$ which restricts to the identity on $\A$. Then there is an induced quasi-isomorphism of $A_\infty$-algebras $\phi^\sigma: \tilde\B^\sigma \rightarrow \B^\sigma$.
\end{lemma}

\proof Our given $\phi$ has components $\phi^{s|1|r}: \A^{\otimes s} \otimes \tilde\B \otimes \A^{\otimes r} \rightarrow \B[-r-s]$. Define
\begin{equation}
\phi^{\sigma,1}(a_+,a_-,b) = (a_+,a_-,\phi^{0|1|0}(b)),
\end{equation}
which is obviously a quasi-isomorphism, and for $d \geq 2$,
\begin{multline}
\phi^{\sigma,d}((a_{d,+},a_{d,-},b_d),\dots,(a_{1,+},a_{1,-},b_1)) \\ =
\Big(0,0, \sum_{i=1}^d \phi^{d-i|1|i-1}(a_{d,+},\dots,a_{i+1,+},b_{i},a_{i-1,-},\dots,a_{1,-})\Big).
\end{multline}
It is straightforward to check the $A_\infty$-homomorphism equations. \qed

\begin{lemma} \label{th:stabilized-trivial}
Suppose that there is an $A_\infty$-bimodule map $\xi: \B/\A \rightarrow \B$ whose composition with the projection $\pi: \B \rightarrow \B/\A$ is a quasi-isomorphism from $\B/\A$ to itself. Then $\B^\sigma$ is quasi-isomorphic to the trivial extension $A_\infty$-algebra $\A \oplus (\B/\A)[-1]$.
\end{lemma}

\proof Consider the $\A$-bimodule $\tilde\B = \A \oplus (\B/\A)$. The inclusion $\iota$ and the map $\xi$ combine to yield an $\A$-bimodule map $\iota \oplus \xi: \tilde\B \rightarrow \B$. Just seen as a map of chain complexes, this sits in a commutative diagram
\begin{equation}
\xymatrix{
0 \ar[r] & \A \ar[r] \ar[d]^{\mathrm{id}} & \tilde\B \ar[r] \ar[d]^{\iota \oplus \xi} & \B/\A \ar[r] \ar[d]^{(\pi \circ \xi)^{0|1|0}} & 0 \\
0 \ar[r] & \A \ar[r] & \B \ar[r] & \B/\A \ar[r] & 0,
}
\end{equation}
hence is itself a quasi-isomorphism. Now consider $\tilde\B$ as an $A_\infty$-algebra, namely the trivial extension of $\A$ by $\B/\A$. By Lemma \ref{th:stabilized-natural}, our bimodule map induces a quasi-isomorphism between the suspensions $\tilde{\B}^\sigma$ and $\B^\sigma$. But Lemma \ref{th:trivial-extension} says that the first of these is quasi-isomorphic to $\A \oplus (\B/\A)[-1]$, which implies the desired result. \qed

\begin{lemma} \label{th:stabilized-split}
Take any $(\A,\B)$. Then there is a map $\xi^\sigma: \B^\sigma/\A^\sigma \rightarrow \B^\sigma$ of bimodules over $\A^\sigma = \A$ whose composition with the projection $\B^\sigma \rightarrow \B^\sigma/\A^\sigma$ is a quasi-isomorphism from $\B^\sigma/\A^\sigma$ to itself.
\end{lemma}

\proof Consider the subspace of $\B^\sigma$ consisting of elements of the form $(a_+,0,b)$. This is a sub-bimodule over $\A^\sigma$. Projection from that bimodule to $\B^\sigma/\A^\sigma$ is an isomorphism, hence admits a unique strict inverse. \qed

The principal consequence of the discussion above is the following:

\begin{proposition} \label{th:double-stabilization}
For any $(\A,\B)$, the double suspension $\B^{\sigma\sigma}$ is an $A_\infty$-algebra quasi-isomorphic to the trivial extension $\A \oplus (\B/\A)[-2]$.
\end{proposition}

\proof Apply Lemma \ref{th:stabilized-split} to $(\A,\B)$, and then Lemma \ref{th:stabilized-trivial} to $(\A^\sigma,\B^\sigma)$. The consequence is that $\B^{\sigma\sigma}$ is quasi-isomorphic to the trivial extension of $\A^\sigma \iso \A$ by the bimodule $(\B^\sigma/\A^\sigma)[-1]$. Projection $\B^\sigma/\A^\sigma = \A \oplus \B[-1] \rightarrow (\B/\A)[-1]$ is a quasi-isomorphism of $\A$-bimodules, and induces a quasi-isomorphism of the associated trivial extension $A_\infty$-algebras. \qed

\section{Fukaya categories\label{sec:fuk}}

Let $(M,\omega)$ be a $2n$-dimensional symplectic manifold, and $(V_1,\dots,V_m)$ an ordered collection of Lagrangian spheres in $M$. To apply Floer-Fukaya theory, in its simplest form, we want to impose a number of restrictions:
\begin{itemize} \itemsep1em
\item $\omega$ should be exact, in fact we want to choose a particular primitive $\omega = d\theta$. To deal with the resulting inevitable non-compactness of $M$, we assume that it can be exhausted by a sequence of relatively compact open subsets $U_1 \subset U_2 \subset \cdots$ with smooth boundaries, such that the Liouville vector field dual to $\theta$ points strictly outwards along each $\partial U_i$.

\item $c_1(M)$ should be trivial, in fact we want to choose a trivialization of the canonical bundle.

\item The $V_i$ should be exact, and we equip them with gradings as well as {\it Spin} structures. Actually, exactness ($\theta|V_i = dF_i$ for some function) is nontrivial only in the lowest dimension $n = 1$. The existence of gradings requires zero Maslov index, which is again nontrivial only for $n = 1$. Finally, the {\it Spin} structure is unique up to isomorphism unless $n = 1$, in which case we always choose the nontrivial {\it Spin} structure on the circle.
\end{itemize}

In this framework, we have well-defined Floer cohomology groups $\mathit{HF}^*(V_i,V_j)$, which are graded vector spaces over our fixed coefficient field $\K$. Following Fukaya (see for instance \cite{seidel04} for an exposition) one constructs an $A_\infty$-category $\B$, with objects $(V_1,\dots,V_m)$, whose morphism spaces are the Floer cochain groups $\mathit{CF}^*(V_i,V_j)$. We sometimes find it convenient to adopt an equivalent point of view. Namely, consider the semisimple ring
\begin{equation}
R = \K e_1 \oplus \cdots \oplus \K e_m, \text{ where $e_i^2 = e_i$, $e_ie_j = 0$ for all $i \neq j$.}
\end{equation}
$A_\infty$-categories with $m$ numbered objects can also be viewed as $A_\infty$-algebras over $R$. Concretely, consider the direct sum
\begin{equation} \label{eq:b}
\B = \bigoplus_{i,j} \mathit{CF}^*(V_i,V_j),
\end{equation}
with its natural $R$-bimodule structure (left multiplication with $e_k$ projects to the direct summand $j = k$, and right multiplication to the direct summand $i = k$). This carries $A_\infty$-multiplications, obtained in the straightforward way from those of the previously mentioned $A_\infty$-category, which are compatible with the $R$-bimodule structure. We will mostly use the point of view of \eqref{eq:b}, but switch back freely to $A_\infty$-categorical language whenever that simplifies the discussion.

As explained in \cite[Theorem 3.2.1.1]{lefevre} or \cite[Lemma 2.1]{seidel04}, we can modify the $A_\infty$-structure on $\B$ to a quasi-isomorphic one which is strictly unital. Assume from now on that this has been done, without changing notation. Concretely, strict unitality means that there is a unique inclusion $R \hookrightarrow \B$, compatible with the $R$-bimodule structure, such that
\begin{equation}
\begin{aligned}
& \mu_\B^1(e_i) = 0, \\
& \mu_\B^2(a,e_i) = a e_i, \\
& \mu_\B^2(e_i,a) = (-1)^{\|a\|-1} e_i a, \\
& \mu_\B^d(\dots,e_i,\dots) = 0 \text{ for all $d \geq 3$.}
\end{aligned}
\end{equation}
In particular, we can then define the directed $A_\infty$-subalgebra of $\B$ to be the subspace
\begin{equation}
\A = R \oplus \bigoplus_{i<j} \mathit{CF}^*(V_i,V_j) \subset \B,
\end{equation}
which is automatically closed under all the $A_\infty$-multiplications. We will need one general property of the pair $(\A,\B)$.

\begin{proposition} \label{th:weak-duality}
The quotient $\B/\A$, with its natural $\A$-bimodule structure, is quasi-isomorphic to $\A^\vee[-n]$, which is the dual of the diagonal bimodule shifted upwards by $n$.
\end{proposition}

\proof We first explain the analogous but much simpler argument on the level of cohomology. There are Poincar{\'e} duality isomorphisms
\begin{equation} \label{eq:floer-poincare}
\mathit{HF}^*(V_i,V_j) \iso \mathit{HF}^{n-*}(V_j,V_i)^\vee
\end{equation}
(generally, the sign of these depends on orientations of the Lagrangian spheres, but in our case there are preferred orientations determined by the gradings). The direct sum of these yields an isomorphism $H(\B) \iso H(\B)^\vee[-n]$, which is compatible with the structure of both as bimodules over $H(\A)$. By definition of the directed subalgebra, the composition
\begin{equation}
H(\A) \hookrightarrow H(\B) \iso H(\B)^\vee[-n] \twoheadrightarrow H(\A)^\vee[-n]
\end{equation}
is zero, which means that we get an induced map $H(\B)/H(\A) = H(\B/A) \rightarrow H(\A)^\vee[-n]$. It is straightforward to verify that this is an isomorphism. On the cochain level, one proceeds as follows:

{\bf Step 1.} {\em There is a quasi-isomorphism of $\A$-bimodules, $\eta: \B \rightarrow \B^\vee[-n]$, which on cohomology induces the maps \eqref{eq:floer-poincare}.}

This is a weak version of the Calabi-Yau property of Fukaya categories. For coefficient fields $\K$ containing $\R$, the strong form of that property (cyclicity) has been proved by Fukaya \cite{fukaya09}. This would be enough for the applications to mirror symmetry, where $\K = \C$. There is actually another argument leading to the weaker form, but which works for arbitrary $\K$; we'll mention it briefly, omitting all details. The general algebraic framework for homomorphism $\B \rightarrow \B^\vee[-n]$ was set up in \cite[Section 5]{tradler01}. A look at the pictures there makes it clear how such a map can be constructed in the Fukaya-theoretic context by looking at moduli spaces of infinite strips with added marked boundary points (strictly speaking, the resulting bimodule map applies to the $\B$ obtained directly from geometry, but of course it carries over to any quasi-isomorphic $A_\infty$-structure).

{\bf Step 2.} {\em The composition of $\A$-bimodule maps
\begin{equation} \label{eq:i-phi-p}
\A \hookrightarrow \B \stackrel{\eta}{\longrightarrow} \B^\vee[-n] \twoheadrightarrow \A^\vee[-n]
\end{equation}
is nullhomotopic}.

Recall that $\A$ is an $A_\infty$-algebra over $R$. All our bimodules and maps between them are compatible with that structure (as in \cite{seidel08}, for instance). Secondly, $\A$ is strictly unital, and therefore every homomorphism between strictly unital $\A$-bimodules is chain homotopic to a strictly unital one. We apply this to $\eta$ and, without changing notation, assume that it is strictly unital. Concretely, this means that
\begin{equation}
\eta^{s|1|r}: \A^{\otimes_R s} \otimes_R \B \otimes_R \A^{\otimes_R r} \rightarrow \B^\vee[-n]
\end{equation}
vanishes whenever any of the $\A$ entries is equal to some $e_i$. Because of this, the $R$-bimodule structure, and directedness, the only possible nonzero component of \eqref{eq:i-phi-p} is the linear part
\begin{equation}
R = \bigoplus_i e_i\A e_i \longrightarrow (\bigoplus_i e_i\A e_i)^\vee = R^\vee[-n]
\end{equation}
which vanishes for degree reasons.

{\bf Step 3.} {\em There is an $\A$-bimodule map $\tilde\eta: \B/\A \rightarrow \A^\vee[-n]$ which fits into a diagram, commutative up to homotopy,
\begin{equation}
\xymatrix{
\B/\A \ar[drr]^{\tilde\eta} && \\
\ar@{->>}[u]
\B \ar[r]^-{\eta} & \B^\vee[-n] \ar@{->>}[r] & \A^\vee[-n]
}
\end{equation}
}

This is a consequence of the previous step, and the fact that a short exact sequence of $A_\infty$-bimodules induces a long exact sequence of morphism spaces. Namely, denoting by ${\mathfrak C}$ the differential graded category of bimodules, we have a long exact sequence
\begin{multline}
\cdots \rightarrow H(\mathit{hom}_{\mathfrak C}(\B/\A,\A^\vee)) \longrightarrow H(\mathit{hom}_{\mathfrak C}(\B,\A^\vee)) \\ \longrightarrow H(\mathit{hom}_{\mathfrak C}(\A,\A^\vee)) \longrightarrow \cdots
\end{multline}

By construction, the cohomology level map induced by $\tilde\eta$ is obtained by taking some splitting of the projection $H(\B) \rightarrow H(\B)/H(\A)$ and composing it with $H(\B) \iso H(\B)^\vee[-n] \rightarrow H(\A)^\vee[-n]$. The result is independent of the choice of splitting, and is an isomorphism, hence $\tilde\eta$ is a quasi-isomorphism. \qed

\section{Lefschetz fibrations\label{sec:p-l}}
Picard-Lefschetz theory can be formulated in various frameworks and degrees of generality. Even though the results here are of a purely symplectic nature, we prefer to keep close to the original algebro-geometric context for expository reasons. Hence, by a Lefschetz fibration we will mean a map
\begin{equation}
p: X \longrightarrow \C
\end{equation}
of the following kind. $X$ is a smooth affine algebraic variety, equipped with a K{\"a}hler form $\omega$ which is exact and complete (in the sense of Riemannian geometry). Additionally, we assume that $X$ carries a holomorphic volume form $\eta$. The function $p$ is regular (which means a polynomial), has only nondegenerate critical points, and additionally satisfies the following Palais-Smale type property:
\begin{equation} \label{eq:palais-smale}
\parbox{30em}{
Any sequence $(x_k)$ of points in $X$ such that $|p(x_k)|$ is bounded and $\|(\nabla p)_{x_k}\| \rightarrow 0$ has a convergent subsequence.
}
\end{equation}
This condition, similar to that of ``tameness'' in singularity theory \cite{broughton88}, implies the well-definedness of symplectic parallel transport maps away from the critical points (because the norm of the parallel transport vector field is essentially the inverse of that of $\nabla p$). In fact, the entire package of symplectic Picard-Lefschetz theory, see for instance \cite[Chapter 3]{seidel04}, applies.

\begin{example} \label{th:laurent}
Take $X = (\C^*)^{n+1}$, with the ``logarithmic'' K{\"a}hler form $\omega = \sum_k d\mathrm{re}(\log(z_k)) \wedge d\mathrm{im}(\log(z_k))$. Let $p \in \C[z_1^{\pm 1},\dots,z_{n+1}^{\pm 1}]$ be a Laurent polynomial satisfying the following two conditions: its Newton polytope contains the origin in its interior; and it is nondegenerate in the sense of \cite{kouchnirenko76}. Then \eqref{eq:palais-smale} holds (compare \cite[Proposition 3.4]{broughton88}, which is an analogue for ordinary polynomials). In fact, these conditions imply the stronger property that $\|\nabla p\|$ is a proper function on $X$.
\end{example}
%rescaling argument

The suspension of a Lefschetz fibration is defined as
\begin{equation}
p^\sigma(x,y) = p(x) - y^2: X^\sigma = X \times \C \longrightarrow \C.
\end{equation}
To be precise, $X^\sigma$ carries the standard product K{\"a}hler form $\omega^\sigma = \omega + \frac{i}{2} \mathit{dy} \wedge \mathit{d\bar{y}}$, and the holomorphic volume form $\eta^\sigma = \eta \wedge dy$. If $p$ satisfies \eqref{eq:palais-smale}, then so does its suspension.

We will now analyze the geometry of the suspension (this material is also covered in \cite[Section 18]{seidel04}, in a closely related setup and with more details). Denote by $X_z$ and $X_z^\sigma$ the fibres of $p$ and $p^\sigma$, respectively, over $z \in \C$. Suppose, just for simplicity of notation, that each fibre of $p$ contains only one critical point, and also that there are no critical values on the positive real half-axis $\R^+ \subset \C$. Projection to the $y$-variable is itself a Lefschetz fibration
\begin{equation} \label{eq:tilde-p}
\tilde{p}: X_0^\sigma \longrightarrow \C.
\end{equation}
Its fibres are $\tilde{p}^{-1}(z) = X_{z^2}$, and the symplectic connection is the pullback of that on $p$ by the double covering of the base. Suppose that we have chosen a distinguished basis $(\gamma_1,\dots,\gamma_m)$ of vanishing paths for $p$, all starting at $z = 0$ and otherwise avoiding $\R^+$. Let $(V_1,\dots,V_m)$ be the associated collection of vanishing cycles, which are Lagrangian spheres in $X_0$. Then
\begin{equation}
(\tilde{\gamma}_1,\dots,\tilde{\gamma}_{2m}) =  (\sqrt{\gamma_1},\dots,\sqrt{\gamma_m},-\sqrt{\gamma_1},\dots,-\sqrt{\gamma_m})
\end{equation}
is a distinguished basis of vanishing paths for $\tilde{p}$, and by our previous observation concerning the symplectic connection, the associated basis of vanishing cycles consists of two copies of that for $p$:
\begin{equation} \label{eq:double}
(\tilde{V}_1,\dots,\tilde{V}_{2m}) = (V_1,\dots,V_m,V_1,\dots,V_m).
\end{equation}
Since there are two instances of each vanishing cycle $V_i$ in that collection, one can join the respective critical points by a matching path $\mu_i$ \cite[Section 16g]{seidel04}, which is simply the composition of $\tilde{\gamma}_i$ and $\tilde{\gamma}_{i+m}$, or equivalently the preimage of $\gamma_i$ under the double branched cover. The matching cycles for these paths form a collection $(V_1^\sigma,\dots,V_m^\sigma)$ of Lagrangian spheres in the total space of $\tilde{p}$, which of course also happens to be the fibre $X^\sigma_0$. These spheres have two interpretations. On one hand, each $V_i$ bounds a Lefschetz thimble, which is a Lagrangian ball in $X$, fibered over $\gamma_i$. The preimage of this ball under the branched covering $X_0^\sigma \rightarrow X$ is $V_i^\sigma$ (this follows directly from the definition of matching cycle). Alternatively, we can use $\gamma_i$ as a vanishing path for $p^\sigma$ itself, and then $V_i^\sigma$ is the vanishing cycle associated to that path (this is a standard braid monodromy argument, compare \cite[Lemma 18.2]{seidel04}).

Passing to Floer-Fukaya theory, we now have the following algebraic structures:
\begin{itemize} \itemsep1em
\item $\B$ is the $A_\infty$-algebra over $R = \K^m$ associated to $V_1,\dots,V_m \subset X_0$. Here, we assume that the vanishing cycles have been equipped with gradings and {\em Spin} structures. As in Section \ref{sec:fuk}, we also assume that $\B$ has been made strictly unital, so that we have a directed $A_\infty$-subalgebra $\A$.

\item Likewise, let $\tilde\B$ be the $A_\infty$-algebra over $\tilde{R} = \K^{2m}$ associated to $\tilde{V}_1,\dots,\tilde{V}_{2m} \subset X_0$, where the gradings are borrowed from those of the $V_i$, and $\tilde\A \subset \tilde\B$ its directed subalgebra.
\end{itemize}
The algebraic relation between these two is straightforward, in view of \eqref{eq:double}. In categorical terms, $\tilde\B$ is obtained from $\B$ by introducing two isomorphic copies of each existing object. Equivalently, as an algebra it is given by
\begin{equation}
\tilde\B = \B \otimes_\K \mathit{mat}_2(\K),
\end{equation}
where $\mathit{mat}_2(\K)$ is the associative algebra of matrices of size two. If one thinks in terms of such matrices, then the directed $A_\infty$-subalgebras are related by
\begin{equation} \label{eq:submatrix}
\tilde{\A} = \begin{pmatrix} \A & 0 \\ \B & \A \end{pmatrix} \subset \tilde{\B}.
\end{equation}
Let's give another and even more explicit description. The morphism spaces in $\tilde{\A}$ are
\begin{equation}
e_j \tilde{\A} e_i = \begin{cases}
e_j \A e_i & \text{if $i \leq j \leq m$ or $m < i \leq j$}, \\
e_{j-m} \B e_i & \text{if $i \leq m < j$,} \\
0 & \text{if $i>j$;}
\end{cases}
\end{equation}
and the $A_\infty$-structure is induced from that of $\B$ (more precisely, only the $A_\infty$-structure of $\A$ and the structure of $\B$ as an $A_\infty$-bimodule over $\A$ enter into the construction of $\tilde{\A}$). We introduce two more $A_\infty$-structures:
\begin{itemize} \itemsep1em
\item Think of $\tilde\A$ as an $A_\infty$-category with objects $(\tilde{V}_1,\dots,\tilde{V}_{2m})$, and let $\mathit{Tw}(\tilde\A)$ be the associated category of twisted complexes. Consider the particular twisted complexes
\begin{equation} \label{eq:c-i}
S_i = \mathit{Cone}(e_i: \tilde{V}_i \rightarrow \tilde{V}_{i+m}),
\end{equation}
where the arrow is $e_i \in \mathit{hom}_{\tilde\A}(\tilde{V}_i,\tilde{V}_{i+m}) = \mathit{hom}_\B(V_i,V_i) = e_i \B e_i$, in other words is derived from the identity elements in $\B$. We then define $\B^\sigma_{\mathrm{alg}}$ to be the full $A_\infty$-subcategory of $\mathit{Tw}(\tilde\A)$ with objects $(S_1,\dots,S_m)$. This can also be seen as an $A_\infty$-algebra over $R = \K^m$, as usual, and we denote by $\A^\sigma_{\mathrm{alg}}$ its directed subalgebra.
\item $\B^\sigma_{\mathrm{geom}}$ is the $A_\infty$-algebra associated to $V_1^\sigma,\dots,V_m^\sigma \subset X_0^\sigma$, and (after modifying that to make it strictly unital, as usual) $\A^\sigma_{\mathrm{geom}}$ is its directed $A_\infty$-subalgebra.
\end{itemize}

\begin{lemma} \label{th:alg}
$(\A^\sigma_{\mathrm{alg}},\B^\sigma_{\mathrm{alg}})$ is the suspension $(\A^\sigma,\B^\sigma)$ of the pair $(\A,\B)$.
\end{lemma}

Here, we are implicitly generalizing the material of Sections \ref{sec:def}--\ref{sec:ext} to the context of $A_\infty$-algebras over $R$, which is entirely unproblematic. With this in mind, the Lemma is merely a reformulation of the description given in Remark \ref{th:twisted}. In fact, if we considered the objects $S_i$ as lying in $\mathit{Tw}(\tilde\B) \iso \mathit{Tw}(\B)$, the associated endomorphism algebra would be exactly $\B \otimes \mathit{hom}_\K(C,C)$ as previously considered in Section \ref{sec:def}. Restricting morphisms to $\tilde\A \subset \tilde\B$ yields the subspace $\B^\sigma \subset \B \otimes \mathit{hom}_\K(C,C)$. Of course, the directed subalgebras are then also the same.

\begin{lemma} \label{th:geom}
Provided that $\mathrm{char}(\K) \neq 2$, and that the gradings of the $V_i^\sigma$ have been chosen appropriately, $\B^\sigma_{\mathrm{alg}}$ is quasi-isomorphic to $\B^\sigma_{\mathrm{geom}}$.
\end{lemma}

This is a direct application of results from \cite[Section 10]{seidel04} to the Lefschetz fibration \eqref{eq:tilde-p} (technical tricks used in proving those results are responsible for the restriction on the ground field; there is no reason to believe that this restriction is fundamentally necessary). Figure \ref{fig:2m} shows the basis of vanishing paths for that fibration, as well as the matching paths (in the interest of clarity, the position of the critical points has been distorted somewhat). In such a situation, \cite[Proposition 18.21]{seidel04} says that the full $A_\infty$-subcategory of the Fukaya category $\Fuk(X_0^\sigma)$ consisting of the associated matching cycles is quasi-isomorphically embedded into $\mathit{Tw}(\tilde\A)$. More precisely, the object corresponding to the $i$-th matching cycle $V_i^\sigma$ (with suitably chosen grading) is the mapping cone of the lowest degree nontrivial morphism $\tilde{V}_i \rightarrow \tilde{V}_{i+m}$. But in our case, $H(\mathit{hom}_{\tilde\A}(\tilde{V}_i,\tilde{V}_{i+m})) = HF^*(V_i,V_i) = H^*(S^n;\K)$, so the lowest degree nontrivial morphism is represented by the identity $e_i$, leading to \eqref{eq:c-i}. On general grounds \cite[Theorem 3.2.2.1]{lefevre}, the quasi-isomorphism $\B^{\sigma}_{\mathrm{alg}} \iso \B^{\sigma}_{\mathrm{geom}}$ can be made strictly unital, and then it restricts to a quasi-isomorphism of the associated directed $A_\infty$-subalgebras.
\begin{figure}
\begin{centering}
\begin{picture}(0,0)%
\includegraphics{2m.pstex}%
\end{picture}%
\setlength{\unitlength}{3947sp}%
\begingroup\makeatletter\ifx\SetFigFont\undefined%
\gdef\SetFigFont#1#2#3#4#5{%
  \reset@font\fontsize{#1}{#2pt}%
  \fontfamily{#3}\fontseries{#4}\fontshape{#5}%
  \selectfont}%
\fi\endgroup%
\begin{picture}(3096,1409)(703,-1273)
\put(2126,-736){\makebox(0,0)[lb]{\smash{{\SetFigFont{10}{12.0}{\rmdefault}{\mddefault}{\updefault}{\color[rgb]{0,0,0}$\tilde{\gamma}_{m+1}$}%
}}}}
\put(826,-736){\makebox(0,0)[lb]{\smash{{\SetFigFont{10}{12.0}{\rmdefault}{\mddefault}{\updefault}{\color[rgb]{0,0,0}$\tilde{\gamma}_1$}%
}}}}
\put(2176,-61){\makebox(0,0)[lb]{\smash{{\SetFigFont{10}{12.0}{\rmdefault}{\mddefault}{\updefault}{\color[rgb]{0,0,0}$\mu_1$}%
}}}}
\end{picture}%
\caption{\label{fig:2m}}
\end{centering}
\end{figure}

\begin{theorem} \label{th:suspension}
Assume that $\mathrm{char}(\K) \neq 2$. Let $(\A,\B)$ be the pair of $A_\infty$-algebras associated to $p$, for some choice of vanishing paths. Then the corresponding pair for $p^\sigma$ is quasi-isomorphic to the algebraic suspension $(\A^\sigma,\B^\sigma)$.
\end{theorem}

This follows directly from Lemmas \ref{th:alg} and \ref{th:geom}. Moreover, applying Propositions \ref{th:double-stabilization} and \ref{th:weak-duality}, we obtain the following consequence:

\begin{corollary} \label{th:double}
Let $p^{\sigma\sigma}$ be the double suspension of $p$. Take any basis of vanishing cycles, and let $(\A^{\sigma\sigma},\B^{\sigma\sigma})$ be the associated pair of $A_\infty$-algebras. Suppose that $\mathrm{char}(\K) \neq 2$. Then $\B^{\sigma\sigma}$ is the trivial extension of $\A$ by the bimodule $\A^\vee[-n-2]$ ($n+2$ is the complex dimension of the fibre of $p^{\sigma\sigma}$).
\end{corollary}

\begin{remark}
$p^{\sigma\sigma}: X^{\sigma\sigma} = \C^2 \times X \rightarrow \C$ carries a fibrewise circle action. It seems plausible to think that an appropriate $S^1$-equivariant version of $\B^{\sigma\sigma}$ would recover the information lost during the suspension process (the corresponding statement is true in the elementary topological counterpart of this process, discussed in Section \ref{sec:topology}, as one can see by applying the localization theorem in equivariant cohomology).
\end{remark}

\section{Local mirror symmetry\label{sec:mirror}}

Let $Y$ be any smooth complex projective variety, and $\iota: Y \rightarrow K_Y$ its embedding into the canonical bundle as the zero-section. Take any collection of objects $E_1,\dots,E_m \in D^b(\Coh(Y))$. Let $\CC$ be the $A_\infty$-algebra over $R = \C^m$, defined using suitable resolutions (Dolbeault, Cech, or injective), which underlies
\begin{equation}
H(\CC) = \bigoplus_{i,j=1}^m \mathit{Hom}^*_Y(E_i,E_j),
\end{equation}
where the $\mathit{Hom}$s are morphisms of any degree in the derived category. There is a similar algebra $\D$ for the objects $\iota_* E_i \in D^b_Y(\Coh(K_Y))$.

\begin{theorem}[\protect{Segal \cite[Theorem 4.2]{segal07}, Ballard \cite[Proposition 4.14]{ballard08}}] \label{th:b-side}
$\D$ is quasi-isomorphic to the trivial extension $\CC \oplus \CC^\vee[-\mathrm{dim}_\C(K_Y)]$.
\end{theorem}

Now suppose that the $E_i$ generate $D^b(\Coh(Y))$, which implies that the $\iota_*E_i$ generate $D^b_Y(\Coh(K_Y))$. Consider $\D$ as an $A_\infty$-category with objects $(\iota_*E_1,\dots,\iota_*E_m)$, and let $D(\D)$ be its derived category, defined using twisted complexes as $H^0(\mathit{Tw}(\mathcal{D}))$. In that case
\begin{equation} \label{eq:generate}
D^b_Y(\Coh(K_Y)) \iso D(\mathcal{D}).
\end{equation}

\proof[Proof of Theorem \ref{th:local-hms}] Take a toric del Pezzo surface $Y$ and its mirror $p: X \rightarrow \C$. The mirror is a Laurent polynomial which, for generic choice of coefficients, has nondegenerate critical points and satisfies the conditions from Example \ref{th:laurent} (the specific choice of coefficients is irrelevant, since the space of allowed choices is connected). Fix a basis of vanishing paths for $p$. Let $(V_1,\dots,V_m)$ be the resulting vanishing cycles in some regular fibre $X_z$, and $\A$ the associated directed $A_\infty$-algebra. As mentioned in the Introduction, we know from \cite{auroux-katzarkov-orlov04, auroux-katzarkov-orlov05, ueda06} that there is a full exceptional collection $(E_1,\dots,E_m)$ in $D^b(\Coh(Y))$ such that $\A$ is quasi-isomorphic to the associated $A_\infty$-algebra $\CC$ (in fact, in all these cases one can choose the vanishing cycles so that the $E_i$ form a strong exceptional collection, which simplifies computations somewhat; but that is not strictly relevant for our purpose).

On the other hand, let $L_i = V_i^{\sigma\sigma} \subset H = X^{\sigma\sigma}_z$ be the corresponding vanishing cycles for the double suspension. By Corollary \ref{th:double}, the full subcategory of $\Fuk(H)$ with objects $(L_1,\dots,L_m)$ is quasi-isomorphic to the trivial extension $\A \oplus \A^\vee[-3]$. On the other hand, Theorem \ref{th:b-side} says that $\D$ is quasi-isomorphic to $\CC \oplus \CC^\vee[-3]$. Using \eqref{eq:generate}, we get a full embedding
\begin{equation}
D^b_Y(\Coh(K_Y)) \iso D(\D) \iso D(\A \oplus \A^\vee[-3]) \hookrightarrow
D(\Fuk(H)).
\end{equation}

%\bibliographystyle{plain}
%\bibliography{../../../bib/all,../../../bib/new}
\end{document}